\documentstyle [12pt]{article}
\begin{document}
\title{
Dynamical  $q$-deformation in quantum theory and the
stochastic limit
}
\author{L. Accardi\thanks{
Centro Vito Volterra Universita di Roma Tor Vergata,
Via di Tor Vergata snc-00133 Roma, Italy,
{\tt e-mail: accardi@volterra.mat.utovrm.it}},
~S.V. Kozyrev\thanks{
Institute  of Chemical Physics,
Kosygin Street 4,  117334, Moscow, Russia,
{\tt e-mail: kozyrev@mi.ras.ru}}
~and I.V. Volovich\thanks{
Steklov Mathematical
Institute, Gubkin Street 8, GSP-1, 117966, Moscow, Russia,
{\tt e-mail: volovich@mi.ras.ru}}
}

\maketitle

\centerline{\it Centro Vito Volterra Universita di Roma Tor Vergata}

\begin{abstract}{
A model of particle interacting with quantum field
is considered. The model includes as particular cases the polaron model
and non-relativistic quantum electrodynamics. We show that
the field operators obey $q$-commutation relations with $q$ depending on
time. After the stochastic (or van Hove)  limit, due to the nonlinearity, the
atomic and field degrees of freedom become {\it entangled}  in the sense that
the field and the atomic variables no longer commute but give rise to a new
algebra with new commutation relations replacing the Boson ones. This new
algebra allows to give a simple proof of the fact that the non crossing
half-planar diagrams give the dominating contribution in a weak coupling
regime and to calculate explicitly the correlations associated to the new
algebra.  The above results depend crucially on the fact that  we do not
introduce any dipole or multipole approximation.
}
\end{abstract}

\font\ftitle=cmbx10 scaled\magstep1
\font\smc=cmcsc10
\def\ov{\overline}
\def\la{\langle}
\def\ra{\rangle}
\def\R{\Bbb R}
\def\var{\varepsilon}
\def\om{\omega}
\def\Om{\Omega}
\def\pa{\partial}
\def\C{\Bbb C}
\def\N{\Bbb N}
\def\na{{\nabla}}
\def\tr{{\triangle}}

\section{Introduction}

In recent years it has been great interest to $q$-deformed commutational
relations, see for example \cite{CoonYuBa}-\cite{AreVo}.
In many works $q$-deformed relations are considered as an {\it ad hoc}
deformation of the ordinary commutation relations or as a hidden
symmetry algebra.

In this work we show that the so called collective operators $a_\lambda(t,k)$
in a model of particle interacting with quantum field
satisfy the  $q$-deformed commutation relations (see
(\ref{2.5}), (\ref{2.6}), (\ref{2.4}) below)
where the parameter $q$ depends on time.
The collective operators are natural objects in the
stochastic (van Hove) limit  of the model describing interaction of particle
with quantum field.
The stochastic limit is used to derive  the long time behavior of the system
interacting with reservoir, in particular to derive the master equation
\cite{AcLu}-\cite{AcLuVo97c}.
The main result of this work is that in the stochastic limit the
$q$-deformed commutation relations give rise to the generalized
quantum Boltzmann commutational relations.

We investigate the model describing interaction of non-relativistic particle
with  quantum field.
This model is widely studied in elementary particle physics,
solid state physics, quantum optics, see for example
\cite{Bog}-\cite{Feynman}.
We consider the simplest case in which matter is represented by a
single particle, say an electron, whose position and momentum we denote
respectively by  $q=(q_{1}, q_{2}, q_{3})$ and
$p=(p_{1}, p_{2}, p_{3})$ and satisfy the commutation relations
$ [q_{j}, p_{n} ] = i \delta_{jn}$. The electromagnetic field is described by
Boson operators
$a(k)=(a_{1}(k), a_{2}(k), a_{3}(k));a^{\dag}(k)=
(a^{\dag}_{1}(k), \ldots ,a^{\dag}_3(k))$
satisfying the  canonical commutation relations
$[a_{j}(k),a_{n}^{\dag}(k')]=\delta_{jn}\delta(k-k')$.
The Hamiltonian of a free non relativistic atom interacting with a quantum
electromagnetic field is
\begin{equation}\label{1.1}
H=H_0+\lambda H_I=\int\omega(k)a^{\dag}(k)a(k)dk+{\frac{1}{2}}\,p^2+\lambda
H_I
\end{equation}
where $\lambda$ is a small constant, $\omega (k)=|k|$ and
\begin{equation}\label{1.2}
H_I=\int d^3k(g(k)
p\cdot a^{\dag}(k)  e^{ikq} +
\overline g(k)p\cdot a(k)e^{-ikq}  ) + h.c.
\end{equation}
Here $p\cdot a(k)=\sum_{j=1}^{3}p_j a_j(k)$, $p^2=\sum_{j=1}^{3}p_j^2$,
$a^{\dag}(k)a(k)=\sum_{j=1}^{3}a^{\dag}_j(k)a_j(k)$, $kq=\sum_{j=1}^{3}k_j
q_j$.

The general idea of the stochastic limit is to make the time rescaling
$t\to t/\lambda^2$ in the solution of the Schr\"o\-din\-ger equation in
interaction picture $U^{( \lambda )}_t=e^{itH_0} e^{-itH}$,
associated to the Hamiltonian $H$, i.e.
\begin{equation}\label{1.3}
{ \frac{\partial}{ \partial t} }  U^{( \lambda )}_t
=-i { \lambda  } H_I (t) \ U_t^{
(\lambda )}\qquad , \ U_0^{(\lambda )} = 1
\end{equation}
with $H_I(t)=e^{it H_0}H_Ie^{-itH_0}$ (the {\it evolved interaction
Hamiltonian}).
This leads to the rescaled equation
\begin{equation}\label{1.4}
{ \frac{\partial}{  \partial t} }  U^{( \lambda )}_{t/\lambda^2}
=- {\frac{i}{\lambda}} H_I (t/\lambda^2) \
U_{t/\lambda^2}^{(\lambda )}
\end{equation}
and one wants to study the limits, in a topology to be specified,
\begin{equation}\label{1.5}
\lim_{\lambda\to0}U^{(\lambda)}_{t/\lambda^2}= U_t
\end{equation}
\begin{equation}\label{1.6}
\lim_{\lambda\to0}{\frac{1}{\lambda}}\,H_I\left({\frac{t}{\lambda^2}}\right)
=H_t = \int d^3k \left(g(k)p\cdot b^{\dag}(t,k) +
\overline g(k)p\cdot b(t,k)+h.c.\right)
\end{equation}
Moreover one wants to prove that $U_t$ is the solution of the equation
\begin{equation}\label{1.7}
\partial_t U_t\,=-iH_tU_t\quad;\qquad U_0=1
\end{equation}
The interest of this limit equation is in the fact that many problems
become explicitly integrable.
The stochastic limit of the model (\ref{1.1})-(\ref{1.2})
has been considered in \cite{AcLu}, \cite{AcLuVo97c}, \cite{Gou96},
\cite{Ske96}.

The rescaling $t\to t/\lambda^2$ is equivalent to consider the
simultaneous limit $\lambda\to 0$, $t\to\infty$ under the condition
that $\lambda^2 t$ tends to a constant (interpreted as a new {\it slow scale}
time). This limit captures the main
contributions to the dynamics in a regime, of {\it long times and small
coupling\/} arising  from the cumulative effects, on a large time scale,
of small interactions ($\lambda\to 0$). The physical idea is that,
looked from the slow time scale of the atom, the field looks like a very
chaotic object: a {\it quantum white noise}, i.e. a $\delta$-correlated
(in time) quantum field $b_j^{\dag}(t,k), b_j(t,k)$ also called a {\it master
field}.  If one introduces the dipole approximation the master field is the
usual Boson Fock white noise. Without the dipole approximation the master
field is a completely new type of white noise whose algebra is described by
the relations  \cite{AcLuVo97c}
\begin{equation}\label{1.8}
b_j(t,k)p_n=(p_n+k_n)b_j(t,k)
\end{equation}
\begin{equation}\label{1.9}
b_j(t,k)b_n^{\dag}(t',k')=2\pi\delta(t-t')\delta(\tilde\omega(k)+kp)
\delta(k-k')\delta_{jn}
\end{equation}
\begin{equation}\label{1.10}
\tilde\omega(k):=\omega(k)+{\frac{1}{2}}\,k^2
\end{equation}
Recalling that $p$ is the atomic momentum, we see that
the relation (\ref{1.8})
shows that the atom and the master field are not independent even at a
kinematical level. This is what we call {\it entanglement}.
The relation (\ref{1.9}) is a generalization of the algebra of free
creation--annihilation operators with commutation relations
$$A_iA^{\dag}_j=\delta_{ij}$$
and the corresponding statistics becomes a generalization of the
Boltzmannian (or Free) statistics. This generalization is due to
the fact that the right hand side is
not a scalar but an operator (a function of the atomic momentum). This
means that the relations (\ref{1.8}), (\ref{1.9}) are {\it module commutation relations}.
For any fixed value $\bar p$ of the atomic momentum we get a copy of the
free (or Boltzmannian) algebra. Given the relations
(\ref{1.8}), (\ref{1.9}), (\ref{1.10}), the
statistics of the master field is uniquely determined by the condition
\begin{equation}\label{1.11}
b_j(t,k) \Psi = 0
\end{equation}
where $\Psi $ is the vacuum of the master field, via a module
generalization of the free Wick theorem (this is our Theorem 2 in section (4)
below).

In  Section 2 the dynamically $q$-deformed commutation relations
(\ref{2.5}), (\ref{2.6}), (\ref{2.4}) are
obtained  and the stochastic limit for collective operators is evaluated.
In Section 3 the $n$-point correlation functions of the collective
operators are computed. Finally in Section 4 the stochastic limit
of $n$-point correlation functions  is calculated.

This work was partially supported by INTAS 96-0698
and RFFI-9600312 grants.

\section{Dynamical $q$-deformation }

In order to determine the limit (\ref{1.3}) one rewrites
the rescaled interaction Hamiltonian in terms of some rescaled fields
$a_{\lambda,j}(t,k)$:
\begin{equation}\label{2.1}
{\frac{1}{\lambda}}\,H_I\left({\frac{t}{\lambda^2}}\right)=
\int d^3k p(\overline g(k) a_\lambda(t,k)+  g(k)a^{\dag}_\lambda(t,k)) + h.c.
\end{equation}
where
\begin{equation}\label{2.2}
a_{\lambda,j}(t,k):={\frac{1}{\lambda}}\,e^{i{\frac{t}{\lambda^2}}\,H_0}
e^{-ikq}a_j(k)e^{-i{\frac{t}{\lambda^2}}\,H_0}=
\frac{1}{\lambda}e^{-i{\frac{t}{\lambda^2}}\,
(\tilde\omega(k)+kp)}e^{-ikq}a_j(k)
\end{equation}
Here $\tilde\omega(k)=\omega(k)+{\frac{1}{2}}\,k^2$.
It is now easy to prove
that operators $a_{\lambda,j}(t,k)$ satisfy the following $q$--deformed
module relations,
\begin{equation}\label{2.5}
a_{\lambda,j}(t,k)a^{\dag}_{\lambda,n}(t',k')=
a^{\dag}_{\lambda,n}(t',k')a_{\lambda,j}(t,k)
\cdot q_\lambda(t-t',kk')+
{\frac{1}{\lambda^2}}\,q_\lambda(t-t',\tilde\omega(k)+
kp)\delta(k-k')\delta_{jn}
\end{equation}
\begin{equation}\label{2.6}
a_{\lambda,j}(t,k)p_n=(p_n+k_n)a_{\lambda,j}(t,k)
\end{equation}
where
\begin{equation}\label{2.3}
q_\lambda(t-t',x)=e^{-i{\frac{t-t'}{\lambda^2}}\,x}
\end{equation}
is an oscillating exponent.
This shows that the module $q$--deformation of the commutation relations arise
here as a result of the dynamics and are not put artificially {\it ab
initio}. Now let us suppose that the master field
\begin{equation}\label{2.7}
b_j(t,k)=\lim_{\lambda\to0}a_{\lambda,j}(t,k)
\end{equation}
exist. Then it is natural to conjecture that its algebra shall be
obtained  as the stochastic limit ($\lambda\to0$) of the
algebra  (\ref{2.5}), (\ref{2.6}).
Notice that
the factor $q_\lambda(t-t',x)$ is an oscillating exponent and one easily sees
that
\begin{equation}\label{2.8}
\lim_{\lambda\to0}q_\lambda(t,x)=0\ ,
\qquad\lim_{\lambda\to0}{\frac{1}{\lambda^2}}\,
q_\lambda(t,x)=2\pi\delta(t)\delta(x)
\end{equation}
Thus it is natural
to expect that the limit of (\ref{2.6})
is
\begin{equation}\label{2.9}
b_j(t,k)p_n=(p_n+k_n)b_j(t,k)
\end{equation}
and the limit of (\ref{2.5}) gives the module free relation
\begin{equation}\label{2.10}
b_j(t,k)b_n^{\dag}(t',k')=2\pi\delta(t-t')\delta(\tilde\omega(k)+kp)
\delta(k-k')\delta_{jn}
\end{equation}
Operators $a_{\lambda,j}(t,k)$ also obey the relation
\begin{equation}\label{2.4}
a_{\lambda,j}(t,k)a_{\lambda,n}(t',k')=a_{\lambda,n}(t',k')a_{\lambda,j}(t,k)
q_\lambda^{-1}(t-t',kk')
\end{equation}
In what follows we will not write indexes $j$, $n$ explicitly.
It is clear that the
relation (\ref{2.4}) should disappear after the limit. In fact, if the relation
(\ref{2.4}) would survive in the limit then,
because of (\ref{2.8}), it should give
$b(t,k)b(t',k')=0$, hence also $b^{\dag}(t,k)b^{\dag}(t',k')=0$, so all the
$n$--particle vectors with $n\geq2$ would be zero. But we shall prove
that this is not the case.

An accurate proof of vanishing of relation (\ref{2.4}) looks as follows.
In fact the subject of the stochastic limit is not the algebra
of observables, but the quantum (or algebraic) probablity space.
Quantum probability space is a pair (algebra, state on this algebra).
In the quantum probability space, defined by the algebra
(\ref{2.4}), (\ref{2.5}), (\ref{2.6}) and the vacuum expectation,
we can omit the relation (\ref{2.4}) even before the limit.

Let us explain this fact for simplicity on the example of bosonic algebra.
Consider the algebra ${\cal A}$ with generators $a_i$, $a_j^{\dag}$
and the relations
\begin{equation}\label{2.11}
[a_i, a_j^{\dag}]=\delta_{ij}
\end{equation}
and the state $\langle\cdot\rangle$  on this algebra, equal to the vacuum
expectation in the Fock representation.

It is easy to prove the following lemma.

\noindent{\smc Lemma 1.}\qquad {\sl
In the GNS representation of the algebra ${\cal A}$
with respect to the state $\langle\cdot\rangle$ we have the extra relation
\begin{equation}\label{2.12}
[a_i, a_j]=0
\end{equation}
}

In the language of quantum probability spaces this means that
the quantum probability space (${\cal A}$, $\langle\cdot\rangle$)
is isomorphic to the quantum probability space
(${\cal A'}$, $\langle\cdot\rangle$) where the algebra ${\cal A'}$
is a factor of ${\cal A}$ by the relations (\ref{2.12}).
The isomorphism of quantum probability spaces means the coincidence of all
correlators.

We get, that in the language of quantum probability spaces the relation
(\ref{2.4}) follows from relations (\ref{2.5}), (\ref{2.6})
and the fact, that we use the vacuum expectation.
For this quantum probability space
we investigate the stochastic limit ($\lambda\to 0$).
In the stochastic limit we have to keep the limits of relations
(\ref{2.5}), (\ref{2.6}). The relation (\ref{2.4}) vanishes in the limit,
because the GNS representation of the limiting algebra is realised
not in the symmetric but in the free (or full) Fock space.

To finish the proof we have to prove the existence of the stochastic limit
of $n$--point correlators. This is the subject of the next section.

\section{Calculation of the $n$--point correlator}

In the present section we prove the existence of the limit
of the $q$--deformed  correlators
\begin{equation}\label{3.1}
\langle
a^{\varepsilon_1}_\lambda(t_1,k_1)\dots
a^{\varepsilon_N}_\lambda(t_N,k_N)
\rangle
\end{equation}
where $a^{\varepsilon}$ means either $a$ or $a^{\dag}$ $(\epsilon=0$ for $a$,
$\epsilon=1$ for $a^{\dag}$) and $\langle \cdot \rangle$
denotes vacuum expectation, exist.
Then according to the previous section
the limit of this correlator must be equal to the corresponding
correlator of the master field:
\begin{equation}\label{3.2}
\langle b^{\varepsilon_1}(t_1,k_1)\dots b^{\varepsilon_N}(t_N,k_N)
\rangle
\end{equation}

Let us enumerate annihilators in the product
$
a^{\varepsilon_1}_\lambda(t_1,k_1)\dots
a^{\varepsilon_N}_\lambda(t_N,k_N)
$
as $a_\lambda(t_{m_j},k_{m_j})$, $j=1,\dots J$, and enumerate
creators as  $a^{\dag}_\lambda(t_{m'_j},k_{m'_j})$, $j=1,\dots I$, $I+J=N$.
This means that if $\varepsilon_m=0$ then
$a^{\varepsilon_m}_\lambda(t_m,k_m)=a_\lambda(t_{m_j},k_{m_j})$
for $m=m_j$ (and the analogous condition for $\varepsilon_m=1$).

Let us prove the following lemma.

\medskip

\noindent{\smc Lemma 2}.  {\sl
$$
a_\lambda(t,k)
a^{\varepsilon_1}_\lambda(t_1,k_1)\dots
a^{\varepsilon_N}_\lambda(t_N,k_N)-
$$
$$
-
\prod_{i=1}^I q_{\lambda}^{-1}(t - t_{m_i},k k_{m_i})
\prod_{j=1}^J q_{\lambda}(t - t_{m'_j},k k_{m'_j})
a^{\varepsilon_1}_\lambda(t_1,k_1)\dots
a^{\varepsilon_N}_\lambda(t_N,k_N)
a_\lambda(t,k)=
$$
$$
=\sum_{j=1}^{I}\delta ( k- k_{m'_j})
\frac{1}{\lambda^2}
q_{\lambda}\left( t - t_{m'_j},\tilde\omega (k) +k p \right)
\prod_{m_i<m'_j}
q_{\lambda}\left( t - t_{m'_j},k k_{m_i}\right)
\prod_{m'_i<m'_j}
q^{-1}_{\lambda}\left( t - t_{m'_j},k k_{m'_i}\right)
$$
\begin{equation}\label{recur0}
\prod_{m_i<m'_j}q_{\lambda}^{-1}(t - t_{m_i},k k_{m_i})
\prod_{m'_i<m'_j}q_{\lambda}(t - t_{m'_i},k k_{m'_i})
a^{\varepsilon_1}_\lambda(t_1,k_1)\dots
{\widehat a_\lambda^{\dag}(t_{m'_j},k_{m'_j})} \dots
a^{\varepsilon_N}_\lambda(t_N,k_N)
\end{equation}
Here the notion ${\widehat a_\lambda^{\dag}}$ means that we omit the operator
$a_\lambda^{\dag}$ in this product.
}

{\it Proof}\qquad
The proof of this lemma is by induction over $N$.
The first step of induction is the relation (\ref{2.5}) or (\ref{2.4}).
Given the formula (\ref{recur0})   for $N$, we will prove this formula
for $N+1$.  We consider two cases.

1) The first case: $\varepsilon_{N+1}=0$.  In this case using
(\ref{recur0})  for $N$ and  (\ref{2.4})  we get
$$
a_\lambda(t,k)
a^{\varepsilon_1}_\lambda(t_1,k_1)\dots
a^{\varepsilon_N}_\lambda(t_{N+1},k_{N+1})-
q_{\lambda}^{-1}(t - t_{N+1},k k_{N+1})
$$
$$
\prod_{i=1}^I q_{\lambda}^{-1}(t - t_{m_i},k k_{m_i})
\prod_{j=1}^J q_{\lambda}(t - t_{m'_j},k k_{m'_j})
a^{\varepsilon_1}_\lambda(t_1,k_1)\dots
a_\lambda(t_{N+1},k_{N+1})
a_\lambda(t,k)=
$$
$$
=\sum_{j=1}^{I}\delta ( k- k_{m'_j})
\frac{1}{\lambda^2}
q_{\lambda}\left( t - t_{m'_j},\tilde\omega (k) +k p \right)
\prod_{m_i<m'_j}
q_{\lambda}\left( t - t_{m'_j},k k_{m_i}\right)
\prod_{m'_i<m'_j}
q^{-1}_{\lambda}\left( t - t_{m'_j},k k_{m'_i}\right)
$$
$$
\prod_{m_i<m'_j}q_{\lambda}^{-1}(t - t_{m_i},k k_{m_i})
\prod_{m'_i<m'_j}q_{\lambda}(t - t_{m'_i},k k_{m'_i})
a^{\varepsilon_1}_\lambda(t_1,k_1)\dots
{\widehat a_\lambda^{\dag}(t_{m'_j},k_{m'_j})} \dots
a^{\varepsilon_N}_\lambda(t_N,k_N)
$$
that is exactly (\ref{recur0}) for $N+1$.

2) The second case: $\varepsilon_{N+1}=1$.  In this case using
(\ref{recur0})  for $N$ and (\ref{2.5}) we get
$$
a_\lambda(t,k)
a^{\varepsilon_1}_\lambda(t_1,k_1)\dots
a^{\varepsilon_N}_\lambda(t_{N+1},k_{N+1})-
$$
$$
- \prod_{i=1}^I q_{\lambda}^{-1}(t - t_{m_i},k k_{m_i})
\prod_{j=1}^J q_{\lambda}(t - t_{m'_j},k k_{m'_j})
a^{\varepsilon_1}_\lambda(t_1,k_1)\dots
a^{\varepsilon_N}_\lambda(t_{N},k_{N})
$$
$$
\left(a^{\dag}_\lambda(t_{N+1},k_{N+1})a_\lambda(t,k)
q_{\lambda}(t - t_{N+1},k k_{N+1})+
\delta ( k- k_{N+1})
\frac{1}{\lambda^2}
q_{\lambda}\left( t - t_{N+1},\tilde\omega (k) +k p \right)
\right)=
$$
$$
=\sum_{j=1}^{I}
\delta ( k- k_{m'_j})
\frac{1}{\lambda^2}
q_{\lambda}\left( t - t_{m'_j},\tilde\omega (k) +k p \right)
\prod_{m_i<m'_j}
q_{\lambda}\left( t - t_{m'_j},k k_{m_i}\right)
\prod_{m'_i<m'_j}
q^{-1}_{\lambda}\left( t - t_{m'_j},k k_{m'_i}\right)
$$
$$
\prod_{m_i<m'_j}q_{\lambda}^{-1}(t - t_{m_i},k k_{m_i})
\prod_{m'_i<m'_j}q_{\lambda}(t - t_{m'_i},k k_{m'_i})
$$
$$
a^{\varepsilon_1}_\lambda(t_1,k_1)\dots
{\widehat a_\lambda^{\dag}(t_{m'_j},k_{m'_j})} \dots
a^{\varepsilon_N}_\lambda(t_N,k_N)
a^{\dag}_\lambda(t_{N+1},k_{N+1})
$$
Moving the term
$$
\delta ( k- k_{N+1})
\frac{1}{\lambda^2}
q_{\lambda}\left( t - t_{N+1},\tilde\omega (k) +k p \right)
$$
to the right hand side of this formula
and commuting it with creators and annihilators using (\ref{2.6})
we get (\ref{recur0}) for $N+1$.
This finishes the proof of Lemma 2.

The next theorem describes the form of $N$--point correlator.

\medskip

\noindent{\smc Theorem 1}.  {\sl

i) If the number of creators is not equal to the number of annihilators,
then the  correlator (\ref{3.1})  is equal to zero;

ii) if the number of creators is  equal to the number of annihilators
($N=2n$), then the correlation
function  is equal to the following sum over pair partitions
\begin{equation}\label{3.3}
\sum_{\sigma({\varepsilon})}
\prod^n_{h=1}
\delta ( k_{ m_h}- k_{m'_h})
\frac{1}{\lambda^2} q_{\lambda}\left( \left(t_{ {m}_h} - t_{m'_h}\right),
\left(\tilde\omega (k_{{m_h}}) +k_{m_h} p
+
\sum_{m_\alpha<m_h< {m'}_\alpha }
k_{{ m}_\alpha} \cdot  k_{{ m}_h}
\right)\right)
\end{equation}
$$
\prod_{( {m}_j,m'_j),( {m}_i,m'_i);i,j=1,\dots,n:m_j<m_i<m'_j<m'_i}
q_\lambda(t_{m_i}-t_{m'_j},k_{m_i}\cdot k_{m_j})
$$
where $\sigma(\varepsilon)=\{( {m}_j<m'_j):j=1,\dots,n\}$ is a
partition of $\{1,\dots,2n\}$ associated with
$\varepsilon=(\varepsilon_{1},\dots, \varepsilon_{2n})$.
}

{\it Proof}\qquad
The proof of this theorem is by induction over $n$.
The first step of induction is obvious.
Let us assume the correlator (\ref{3.1}) is expressed by the formula
(\ref{3.3}) for $N=2n-2$ and prove that the same is true for $N=2n$.
We consider $2n$-point correlator
$$
\langle
a^{\varepsilon_1}_\lambda(t_1,k_1)\dots
a^{\varepsilon_{2n}}_\lambda(t_{2n},k_{2n})\rangle
$$
It is easy to see that if this correlator is not equal to zero then
the first operator is annihilator and the last is creator.
Without loss of generality we can consider the case
when the correlator is as follows
\begin{equation}\label{corr1}
\langle
a_\lambda(t_{m_1},k_{m_1})
a^{\varepsilon_2}_\lambda(t_2,k_2)\dots
a^{\varepsilon_{2n-1}}_\lambda(t_{2n-1},k_{2n-1})
a^{\dag}_\lambda(t_{m'_n},k_{m'_n})
\rangle
\end{equation}
From the Lemma 2 follows  the following formula for this correlator
$$
(\ref{corr1})=
\sum_{j=1}^{n}\delta ( k_{ m_1}- k_{m'_j})
\frac{1}{\lambda^2}
q_{\lambda}\left( t_{ {m}_1} - t_{m'_j},\tilde\omega (k_{{m_1}}) +k_{m_1} p
\right)
\prod_{m_i<m'_j<m'_i}
q_{\lambda}\left( t_{ {m}_1} - t_{m'_j},k_{m_1}k_{m_i}\right)
$$
\begin{equation}\label{recur}
\prod_{m_i<m'_j}q_{\lambda}^{-1}(t_{ {m}_1} - t_{m_i},k_{m_1}k_{m_i})
\prod_{m'_i<m'_j}q_{\lambda}(t_{ {m}_1} - t_{m'_i},k_{m_1}k_{m'_i})
\langle {\widehat a_\lambda(t_{m_1},k_{m_1})}\dots
{\widehat a_\lambda^{\dag}(t_{m'_j},k_{m'_j})} \dots \rangle
\end{equation}
The product
$
\prod_{m_i<m'_j<m'_i}
q_{\lambda}\left( t_{ {m}_1} - t_{m'_j},k_{m_1}k_{m_i}\right)
$
in   (\ref{recur}) arise from the products
$$
\prod_{m_i<m'_j}
q_{\lambda}\left( t - t_{m'_j},k k_{m_i}\right)
\prod_{m'_i<m'_j}
q^{-1}_{\lambda}\left( t - t_{m'_j},k k_{m'_i}\right)
$$
in   (\ref{recur0}) due to cancellation of corresponding terms
because of $\delta$-functions  $\delta(k_{m_i}-k_{m'_i})$
in the correlator (\ref{3.3}) for $N=2n-2$.   We have
$$
\prod_{m_i<m'_j}
q_{\lambda}\left( t - t_{m'_j},k k_{m_i}\right)
\prod_{m'_i<m'_j}
q^{-1}_{\lambda}\left( t - t_{m'_j},k k_{m'_i}\right)=
$$
$$
=\prod_{m'_i<m'_j}
q_{\lambda}\left( t - t_{m'_j},k k_{m_i}\right)
\prod_{m_i<m'_j<m'_i}
q_{\lambda}\left( t - t_{m'_j},k k_{m_i}\right)
\prod_{m'_i<m'_j}
q^{-1}_{\lambda}\left( t - t_{m'_j},k k_{m_i}\right)=
$$
$$
=\prod_{m_i<m'_j<m'_i}
q_{\lambda}\left( t - t_{m'_j},k k_{m_i}\right)
$$
Let us prove now that the (\ref{recur}) is equal in fact to (\ref{3.3}).
This will give a proof of the theorem.
We have
$$
\prod_{m_i<m'_j}q_{\lambda}^{-1}(t_{ {m}_1} - t_{m_i},k_{m_1}k_{m_i})
\prod_{m'_i<m'_j}q_{\lambda}(t_{ {m}_1} - t_{m'_i},k_{m_1}k_{m'_i})=
$$
\begin{equation}\label{threeprod}
=\prod_{m'_i<m'_j}q_{\lambda}^{-1}(t_{ {m}_1} - t_{m_i},k_{m_1}k_{m_i})
\prod_{m_i<m'_j<m'_i}q_{\lambda}^{-1}(t_{ {m}_1} - t_{m_i},k_{m_1}k_{m_i})
\prod_{m'_i<m'_j}q_{\lambda}(t_{ {m}_1} - t_{m'_i},k_{m_1}k_{m'_i})
\end{equation}
because $m_i<m'_i$.
From (\ref{3.3}) for $2n-2$ we have
$k_{m_i}=k_{m'_i}$.
By using this and combining the first product with the third
we get
\begin{equation}\label{twoprod}
\prod_{m'_i<m'_j}q_{\lambda}(t_{ {m}_i} - t_{m'_i},k_{m_1}k_{m_i})
\prod_{m_i<m'_j<m_i'}q_{\lambda}^{-1}(t_{ {m}_1} - t_{m_i},k_{m_1}k_{m_i})
\end{equation}
Using the change of variables in the second product in (\ref{twoprod})
$$
t_{ {m}_1} - t_{m_i}=(t_{ {m}_1} - t_{m'_j})-(t_{ {m'}_j} - t_{m_i})
$$
and the property $k_{m_1}=k_{m'_j}$ we get that (\ref{threeprod}) equals
$$
\prod_{m'_i<m'_j}q_{\lambda}(t_{ {m}_i} - t_{m'_i},k_{m_1}k_{m_i})
\prod_{m_i<m'_j<m_i'}q_{\lambda}^{-1}(t_{ {m}_1} - t_{m'_j},k_{m_1}k_{m_i})
\prod_{m_i<m'_j<m_i'}q_{\lambda}(t_{ {m}_i} - t_{m'_j},k_{m'_j}k_{m_i})
$$
Substituting this into the formula (\ref{recur}) we get
\begin{equation}\label{steptwo}
\sum_{j=1}^{n}\delta ( k_{ m_1}- k_{m'_j})
\frac{1}{\lambda^2}
q_{\lambda}\left( t_{ {m}_1} - t_{m'_j},\tilde\omega (k_{{m_1}}) +k_{m_1} p
\right)
\end{equation}
$$
\prod_{m'_i<m'_j}q_{\lambda}(t_{ {m}_i} - t_{m'_i},k_{m_1}k_{m_i})
\prod_{m_i<m'_j<m_i'}q_{\lambda}(t_{ {m}_i} - t_{m'_j},k_{m'_j}k_{m_i})
\langle {\widehat a_\lambda(t_{m_1},k_{m_1})}\dots
{\widehat a_\lambda^{\dag}(t_{m'_j},k_{m'_j})} \dots \rangle
$$
Here the notion $\langle\dots{\widehat a}\dots\rangle$
means that we omit the operator ${\widehat a}$ in this correlation function.
For $2n-2$--point correlator in (\ref{steptwo})
we use the formula (\ref{3.3}) for
$\varepsilon-\{m_1, m'_j \}$:
$$
\langle {\widehat a_\lambda(t_{m_1},k_{m_1})}\dots
{\widehat a_\lambda^{\dag}(t_{m'_j},k_{m'_j})} \dots \rangle=
$$
\begin{equation}\label{step2n-2}
=\sum_{\sigma({\varepsilon-\{m_1, m'_j \}})}
\prod^{n-1}_{h=1}
\delta ( k_{ n_h}- k_{n'_h})
\frac{1}{\lambda^2} q_{\lambda}\left( \left(t_{ {n}_h} - t_{n'_h}\right),
\left(\tilde\omega (k_{{n_h}}) +k_{n_h} p
+
\sum_{n_\alpha <n_h< {n'}_\alpha}
k_{{ n}_\alpha} \cdot  k_{{ n}_h}
\right)\right)
\end{equation}
$$
\prod_{( {n}_j,n'_j),( {n}_i,n'_i);i,j=1,\dots,n-1:n_j<n_i<n'_j<n'_j}
q_\lambda(t_{n_i}-t_{n'_j},k_{n_i}\cdot k_{n'_j})
$$
where $\sigma(\varepsilon-\{m_1, m'_j \})=\{( {n}_j<n'_j):j=1,\dots,n-1\}$
is a partition (without one pair) of $\{1,\dots,2n\}$ associated with
$\varepsilon-\{m_1, m'_j \}$.
The indices $n_h$ correspond to annihilators, $n'_h$ correspond to creators.

Substituting (\ref{step2n-2})  into (\ref{steptwo}) we get
\begin{equation}\label{steplast}
\sum_{j=1}^{n}\delta ( k_{ m_1}- k_{m'_j})
\frac{1}{\lambda^2}
q_{\lambda}\left( t_{ {m}_1} - t_{m'_j},\tilde\omega (k_{{m_1}}) +k_{m_1} p
\right)
\end{equation}
$$
\prod_{m'_i<m'_j}q_{\lambda}(t_{ {m}_i} - t_{m'_i},k_{m_1}k_{m_i})
\prod_{m_i<m'_j<m_i'}q_{\lambda}(t_{ {m}_i} - t_{m'_j},k_{m'_j}k_{m_i})
$$
$$
\sum_{\sigma({\varepsilon-\{m_1, m'_j \}})}
\prod^{n-1}_{h=1}
\delta ( k_{ n_h}- k_{n'_h})
\frac{1}{\lambda^2} q_{\lambda}\left( \left(t_{ {n}_h} - t_{n'_h}\right),
\left(\tilde\omega (k_{{n_h}}) +k_{n_h} p
+
\sum_{n_\alpha <n_h< {n'}_\alpha}
k_{{ n}_\alpha} \cdot  k_{{ n}_h}
\right)\right)
$$
$$
\prod_{( {n}_j,n'_j),( {n}_i,n'_i);i,j=1,\dots,n-1:n_j<n_i<n'_j<n'_j}
q_\lambda(t_{n_i}-t_{n'_j},k_{n_i}\cdot k_{n'_j})
$$
It is easy to see that
\begin{equation}\label{sums}
\sum_{j=1}^n \sum_{\sigma({\varepsilon-\{m_1, m'_j \}})}
=\sum_{\sigma({\varepsilon})}
\end{equation}
Using (\ref{sums}) and
combining the first product in (\ref{steplast}) with the third
and the second product with the fourth we obtain (\ref{3.3}).

This finishes the proof of the theorem.

Let us analize the behavior of the $n$--point correlator in the stochastic
limit. It is easy to see, that the stochastic
limit of pairings
exists and equals to the product of $\delta$-functions of different arguments.
Product of pairing cannot spoil the convergence, because this product
in any case can be considered as the product of two terms. The first term
can result only in shift in the $\delta$-functions, corresponding to pairings.
The second term is an oscillating exponent and vanish in the stochastic limit.

We have proved that the stochastic limit of the $n$--point correlator exists.

\section{The QED module  Wick theorem}

We proved that
$b(t,k)$
satisfy the following  free QED module  algebra relations
\begin{equation}\label{4.1}
b(t,k_1)b^{\dag}(\tau,k_2)=2\pi\delta(t-\tau)\delta(\tilde\omega(k_1)+k_1
p){\delta(k_1-k_2)}
\end{equation}
\begin{equation}\label{4.2}
b(t,k)p=(p+k)b(t,k)
\end{equation}
and  the functional $\langle \cdot \rangle$ is the vacuum expectation.
Let us prove the following module analog of the Wick theorem.

\noindent{\smc Theorem 2.}\qquad {\sl
The limit correlation functions
exist always and

i) if the number of creators is not equal to the number of annihilators,
then the correlator (\ref{3.2})  is equal to zero (even before
the limit);

ii) if the number of creators is  equal to the number of annihilators
($N=2n$), then the limit  is equal to the following
\begin{equation}\label{4.3}
\prod^n_{h=1}
\delta ( k_{ m_h}- k_{m'_h})
2\pi \delta(t_{ {m}_h} - t_{m'_h})\delta
\left(\tilde\omega (k_{{m_h}}) +k_{m_h} p
+
\sum_{m_\alpha <m_h< {m'}_\alpha}
k_{{ m}_\alpha} \cdot  k_{{ m}_h}
\right)
\end{equation}
where $\{( {m}_j<m'_j):j=1,\dots,n\}$ is the unique non-crossing
partition of $\{1,\dots,2n\}$ associated with 
$\varepsilon=(\varepsilon_{1},\dots, \varepsilon_{2n})$.
}

\medskip

{\it Proof}\qquad
The proof is done by computing the correlation functions using the commutation
relations listed above. We investigate the correlator
$$
\langle
b^{\epsilon_{1}}(t_{1},k_{1})\dots
\langle b^{\epsilon_{N}}(t_N,k_N)
\rangle
$$
At first
we simplify this correlator using (\ref{4.1}). Obtained $\delta$-functions
we will move through $b^{\epsilon}(t,k)$, using (\ref{4.2}). We will iterate
this procedure before monomial will take normally ordered form.
Because the functional  $\langle \cdot \rangle$ is equal to vacuum
expectation, only $\delta$-functions will survive.

The pairing
$b(t_{m'_h},k_{m'_h})b^{\dag}(t_{m_h},k_{m_h})$
equal
\begin{equation}\label{4.4}
\delta ( k_{ m'_h}- k_{m_h})
2\pi \delta(t_{ m'_h} - t_{m_h} )
\delta
\left(\tilde\omega (k_{{m_h}}) +k_{m_h} p
\right)
\end{equation}
and the  relation (\ref{4.2}) gives the term
$
\sum_{m_\alpha <m_h< {m'}_\alpha} k_{{ m}_\alpha} \cdot  k_{{ m}_h}
$
in the phase shift
(an argument of the last $\delta$-function in (\ref{4.3})),
arising from moving of this $\delta$-function through $b^{\epsilon}(t,k)$.
This finishes the proof of theorem 2.
The alternative proof can be given by calculation of stochastic limit
of correlator given by the Theorem 1 (for dynamically $q$-deformed algebra).
\bigskip

{\bf Acknowlegements}
\medskip

S.Kozyrev and I.Volovich are grateful to L.Accardi and Centro Vito Volterra
for kind hospitality.

\end{document}